\input amstex
\documentstyle{amsppt}
\input bull-ppt
\keyedby{guest/amh}   

\refstyle{A}

\predefine\preBbb{\Bbb}
\redefine\Bbb{\bold}
\define\C{\Bbb C}
\define\R{\Bbb R}
\define\Z{\Bbb Z}
\predefine\prell{\ll}
\redefine\ll{\lq\lq}
\predefine\prerr{\rr}
\redefine\rr{\rq\rq\ }
\define\rrr{\rq\rq}

\define\Hol{\operatorname {Hol}}
\define\Map{\operatorname {Map}}
\define\Hom{\operatorname {Hom}}
\predefine\premin{\min}
\redefine\min{\operatorname{min}}
\define\SP{\operatorname{Sp}}
\topmatter
\cvol{31}
\cvolyear{1994}
\cmonth{October}
\cyear{1994}
\cvolno{2}
\shorttitle{Configuration spaces on a toric variety}
\cpgs{191-196}
\title Configuration spaces and the space \\
of  rational curves on a toric variety
\endtitle 
\author M. A. Guest
\endauthor
\address Department of Mathematics, University of
Rochester, Rochester, New York 14627\endaddress
\ml 
gues\@db1.cc.rochester.edu\endml
\date November 15, 1993\enddate
\abstract 
The space of holomorphic maps from $S^2$ to a complex 
algebraic 
variety $X$, i.e. the space of parametrized rational 
curves on $X$, 
arises in several areas of geometry. It is a well known 
problem
 to determine an integer 
$n(D)$ such that the inclusion of this space in the 
corresponding 
space of continuous maps induces isomorphisms of homotopy 
groups 
up to dimension $n(D)$, where $D$ denotes the homotopy 
class of the 
maps.   The solution to this problem is known for an 
important but 
special class of varieties, the generalized flag 
manifolds: such an 
integer may be computed, and $n(D)\to\infty$ as 
$D\to\infty$. 
We consider the problem for another class of varieties, 
namely, toric 
varieties. For smooth toric varieties and certain singular 
ones,  $n(D)$ may be 
computed, and $n(D)\to\infty$ as $D\to\infty$. For other 
singular toric
varieties, however, it turns out that $n(D)$ cannot 
always be made arbitrarily large by a suitable choice of 
$D$.\endabstract
\endtopmatter

\document

The space of all holomorphic maps from the Riemann sphere 
$S^2=\C\cup\infty$ to a complex projective variety $X$ 
will be 
denoted by $\Hol(S^2,X)$. (In this context, \ll 
holomorphic\rr means 
the same as \ll algebraic\rr or \ll rational\rr, but we 
use the first 
term so as not to preclude more general situations later.) 
The global 
topological properties of $\Hol(S^2,X)$ are relevant to 
various 
problems in algebraic topology, algebraic geometry, 
differential 
geometry, and mathematical physics. 

When $X=S^2$, $\Hol(S^2,X)$ is just the space of rational 
functions 
--- topologized as a subspace of the space $\Map(S^2,X)$ of 
continuous (or smooth) maps.
It is elementary that each connected component of 
$\Hol(S^2,S^2)$ is 
indexed by a non-negative integer (the degree, or 
the number of poles, of a map). The next topological 
invariant, the 
fundamental group, was computed independently by Epshtein 
(\cite{Ep}) and Jones (see \S 6 of \cite{Se}): if 
$\Hol_d(S^2,S^2)$ 
denotes the maps of degree $d$ and $\Hol^\ast_d(S^2,S^2)$ 
the maps 
of degree $d$ which fix the point $\infty$, then 
$\pi_1\Hol_d(S^2,S^2)$ is cyclic of order $2d$, and 
$\pi_1\Hol^\ast_d(S^2,S^2)\cong\Z$. Now, it happens that 
these agree 
with the fundamental groups of the corresponding spaces 
$\Map_d(S^2,S^2)$, $\Map^\ast_d(S^2,S^2)$. An \ll 
explanation\rr 
for this coincidence is provided by Morse theory, for 
$\Hol_d(S^2,S^2)$ is 
the space of absolute minima of the energy functional 
$E:\Map_d(S^2,S^2)\to\R$, $f\mapsto\int_{S^2}\vert 
df\vert^2$ (see 
\cite{EL}). Thus (cf. \S 22 of \cite{Mi}), one expects 
the space of holomorphic maps 
to have the same homotopy groups as the space of 
continuous maps, up to some 
dimension. (This applies more generally to $\Hol(S^2,X)$ 
for any 
K\"ahler manifold $X$ and is therefore of interest from 
the point of 
view of harmonic maps or minimal immersions of $S^2$ into 
such 
spaces.)  Although Morse theory does not apply naively to 
the energy 
functional here, such a \ll comparison theorem\rr was 
established without appealing to Morse theory by Segal, 
who showed in 
\cite{Se} that the inclusions
$$
\Hol_d(S^2,S^2)\to\Map_d(S^2,S^2)\ \text{\ \ \ and\ \ \ }\ 
\Hol^\ast_d(S^2,S^2)\to\Map^\ast_d(S^2,S^2)
$$
induce isomorphisms of homotopy groups $\pi_i$  for $0\le 
i<d$. 
In particular, these maps also induce isomorphisms in
integral homology groups in the same range. From the point 
of view 
of algebraic topology, this leads to an interesting 
development
(see \cite{CCMM}, and also \cite{CS, Va}), because the 
homology of 
$\Map^\ast_d(S^2,S^2)=\Omega^2_d S^2$ is well known. This 
homology
has a natural algebraic filtration,  and one of the main 
results of \cite{CCMM} is 
that the filtration is realized geometrically by the spaces 
$\Hol^\ast_d(S^2,S^2)$.  It should be mentioned that 
Segal's 
result is also relevant 
to, and was motivated partly by,  a problem in control 
theory (cf. \cite{BD}).

As predicted in \cite{Se}, the comparison theorem has been 
 extended to the 
inclusion  
$\Hol^\ast_d(S^2,X)\to\Map^\ast_d(S^2,X)$ when $X$ is any 
generalized flag 
manifold (i.e. $X=G^{\C}/P$ where $G^{\C}$ is a complex 
semisimple Lie group and 
$P$ is a parabolic subgroup).  This was carried out in 
\cite{Gu1, Ki1, Gr, MM1, MM2,  BHMM2}. It has even 
been extended to the case where $X=\Omega G$ 
(\cite{Gr, Ta, BHMM1, Ki2}; see also \cite{Gu2}), which is 
an 
{\it infinite dimensional} generalized flag manifold. 
Apart from its 
relevance to control theory and to harmonic maps, this 
result is of 
interest from the point of view of Yang-Mills theory, as 
$\Hol^\ast(S^2,G^{\C}/P)$ may be interpreted as a moduli 
space of 
monopoles. The space $\Hol^\ast(S^2,\Omega G)$ may be 
interpreted as 
a moduli space of $G$-instantons over $S^4$, and the 
comparison 
theorem in this case was conjectured by Atiyah and Jones 
\cite{AJ}
in the early  days of mathematical Yang-Mills theory.

With the notable exception of \cite{Se}, the work in 
proving the 
comparison theorem for $X=G^{\C}/P$ employed traditional 
methods of 
algebraic topology and relied either on the results of 
\cite{Se} or on knowing 
beforehand the homology of 
the spaces of continuous maps. It is natural to ask 
whether the theorem 
can be extended to larger classes of complex manifolds 
(or varieties or analytic 
spaces), and it is natural to ask for a more conceptual 
proof. There 
are some obvious limitations; for example, there may be no 
(or only 
finitely many) non-constant holomorphic maps from $S^2$ to 
$X$. In 
the case where $\Hol(S^2,X)$ is large, a conceptual (Morse 
theoretic) 
proof may ultimately be provided by the analytical methods 
of 
Donaldson, Taubes, and Uhlenbeck (cf. \S 4 of \cite{Uh}) 
or by the Gromov-Floer theory developed in 
\cite{CJS}. So far, however,  this goal has proved 
elusive.  

In view of this, it is of interest to investigate concrete 
examples beyond the class of generalized flag manifolds, 
and we shall 
do this by establishing the comparison theorem for a large 
class of 
{\it toric varieties}. Our method is also based on 
\cite{Se}, but it
uses the ideas rather than the results of \cite{Se} and is 
therefore quite
self-contained. Moreover,  our method  represents a 
technical simplification of the 
method of \cite{Se} (which in turn is simpler than the more 
traditional methods referred to above). It also makes 
explicit 
the fundamental idea that the functor 
$\C\mapsto \pi_i\Hol(\C\cup\infty,X)$
behaves like a generalization of a homology  functor.

Before stating our results, we recall briefly 
the definition of a toric variety:
a toric variety is an irreducible normal algebraic variety 
$X$ which 
admits an (algebraic) action of an algebraic torus 
$T^{\C}=\C^\ast\times\dots\times\C^\ast$ such that the 
orbit 
$T^{\C}\cdot \ast$ of some point $\ast\in X$ is a densely 
embedded copy 
of $T^{\C}$. Such varieties were considered in the 1970s 
in connection with 
compactifications in algebraic geometry, 
but they have only recently been 
introduced to wider audiences \cite{Od, Fu}. Their 
attractiveness derives from the fact that they are 
amenable to study 
by combinatorial methods. This is because a toric variety 
$X$ is 
characterized by its \ll fan\rr $\Delta$, a collection of 
convex cones 
in $\R^k$ (where $k=\dim_{\C} T^{\C}$).  In the case of 
projective 
algebraic varieties,  toric varieties are characterized as 
being 
definable by means of \ll monomial equations\rrr.

Our first result gives a description of $\Hol^\ast(S^2,X)$ 
as a \ll 
labelled configuration space\rrr.  (The basepoint 
condition here is taken to be 
$f(\infty)=\ast$.) To establish notation, let $Q(\C;M)$ 
denote the space of configurations of distinct points in 
$\C$ which 
have labels in a partial monoid $M$. This space is 
topologized in such 
a way that two labelled points $z_1,z_2$ are allowed to 
collide if the 
sum of their labels $m_1,m_2$ is defined in $M$ (the 
result of such a 
collision would be a new point with label $m_1+m_2$); 
otherwise the 
collision is prohibited. To a fan $\Delta$ we associate a 
partial 
monoid $M_{\Delta}$, which is a certain subset of 
$\Hom(G,\Z)$, the dual of the 
group $G$ of $T^{\C}$-equivariant Cartier divisors on $X$. 
 There is a 
homomorphism $h:\pi_2 X\to \Hom(G,\Z)$.

\proclaim{Proposition 1} Let $D\in\pi_2 X$. Then the space 
$\Hol_D^\ast(S^2,X)$ is diffeomorphic to the configuration 
space
$$
Q_D(\C;M_{\Delta})=
\{\ \{(z_i,m_i)\}_i\in Q(\C;M_{\Delta})\ \vert\ \sum_i 
m_i=h(D)\}.
$$
\endproclaim

This generalizes the fact that a basepoint-preserving 
rational function is 
determined by its poles and zeros. Similar configuration 
spaces have appeared in algebraic topology since the 
work of Dold and Thom \cite{DT} on symmetric products. 
Indeed, the 
symmetric product $\SP^k(U)$ of a space $U$ is the 
simplest example 
of this type; it may be considered as a subspace of the 
configuration 
space $Q(U;M_k)$, where $M_k$ is the partial monoid 
$\{1,2,\dots,k\}$. The main result of \cite{DT} is that 
the functor 
$U\mapsto \pi_i \SP^k(U)$ agrees with the ordinary 
homology functor $H_i$,
for all $i$  up to some dimension,
on nice spaces $U$. The same method shows that the functor
$U\mapsto \pi_i Q_D(U;M_{\Delta})$ satisfies the axioms 
for a
generalized homology theory, for a range of values of $i$, 
on an
appropriate category.

In the case of a smooth toric variety $X$, $\pi_2 X$ may 
be identified 
naturally with a subgroup of $\Z^u$, where $u$ is the 
number of 
codimension-one orbits of $T^{\C}$ in $X$, so we may write 
$D=(d_1,\dots,d_u)$. The integer $d_i$ represents the 
intersection 
number of the 
map with the closure of the $i$-th orbit. We then have:

\proclaim{Theorem 2} Let $X$ be a smooth toric variety. Then
the inclusion $\Hol^\ast_D(S^2,X)\to 
\Map^\ast_D(S^2,X)$ induces isomorphisms of homotopy groups 
$\pi_i$ for $0\le i< \min\{d_1,\dots$, $d_u\}$.
\endproclaim

The proof has two main ingredients. One is similar to an 
argument of 
\cite{Se}:
it uses the homology-like properties of the functor 
$\C\mapsto\pi_i Q_D(\C\cup\infty;M_\Delta)$ to show that 
one obtains a 
homotopy equivalence in the limit 
$D\to\infty$. The other ingredient is the fact that there 
are \ll stabilization 
maps\rr $Q_{d_1,\dots,d_u}(\C;M_{\Delta})\to 
Q_{d_1^\prime,\dots,d_u^\prime}(\C;M_{\Delta})$ (where 
$d_i^\prime\ge d_i$) which induce isomorphisms of homotopy 
groups 
up to dimension $\min\{d_1,\dots,d_u\}$. Here the method 
differs 
significantly from that of \cite{Se}; it was introduced in 
\cite{GKY}. The idea is to reduce to the corresponding 
statement for symmetric products, which is elementary. 
This has the advantage 
that the 
passage from homology to homotopy is easier than in 
\cite{Se}.

When $X$ is singular, the situation is considerably more 
complicated. 
 Nevertheless, 
we can give a method for determining the extent to which 
the 
comparison theorem holds, based upon the choice of a toric 
resolution 
$\hat X$ of $X$. Such a resolution may be obtained by 
subdividing the 
fan $\Delta$ by \ll inserting rays\rrr. (The precise 
meaning of this is 
explained in \S 2.7 of \cite{Fu}.) The simplest 
application of our method gives:

\proclaim{Theorem 3} Let $X$ be a toric variety.  Assume 
that a toric resolution 
$\hat X$  of $X$ may be found by inserting a single ray in 
the fan $\Delta$ of 
$X$. Then  the inclusion $\Hol^\ast_D(S^2,X)\to 
\Map^\ast_D(S^2,X)$ 
induces isomorphisms of homotopy groups $\pi_i$ for $0\le 
i<n(D)$, 
where $n(D)\to\infty$ as $D\to\infty$.
\endproclaim

In the situation of Theorem 3, we can in fact give 
an explicit formula for 
the integer $n(D)$, but we are not able to give a formula 
in 
the case of a general 
toric variety,
 partly because the integer $n(D)$ depends on the choice of 
resolution (which is far from unique). However, it is 
possible to use 
this method to determine 
$n(D)$ in concrete examples, two of which are given below.

\ex{Example 1} (See \cite{Gu3}.) The subvariety $X$ of $\C 
P^3$ defined by the monomial equation $z_2^2=z_1z_3$ has 
only one 
singular point, which may be resolved by inserting a 
single ray in the 
fan of $X$. Here, $\pi_2 X\cong\Z$, and the method of 
Theorem 3 
shows that the inclusion $\Hol^\ast_d(S^2,X)\to 
\Map^\ast_d(S^2,X)$ 
induces isomorphisms of homotopy groups $\pi_i$ for $0\le 
i<[d/3]$.
\endex

\ex{Example 2} The weighted projective space $X=P(1,2,3)$ 
is by definition the quotient space of $\C^3-\{0\}$ by the 
action 
$u\cdot[z_0;z_1;z_2]=[uz_0;u^2z_1;u^3z_2]$ of $\C^\ast$. 
There
are two singular points, which may be resolved by 
inserting three 
rays. Again, $\pi_2X\cong\Z$, and three applications of the 
above method lead to the conclusion 
that the inclusion 
$\Hol_d^\ast(S^2,X)\to\Map_d^\ast(S^2,X)$ induces
isomorphisms of homotopy groups up to dimension 
$[\tfrac25d]$.
\endex

We do not claim, however, 
 that the method {\it always} leads to an integer $n(D)$ 
such that $n(D)\to\infty$ as $D\to\infty$.
Another source of 
trouble is the choice of basepoint in $X$ --- this is 
relevant 
as toric varieties are not in general homogeneous spaces. 
So far we 
have assumed that the basepoint of $X$ is a point $\ast$ 
such that the 
closure of $T^{\C}\cdot\ast$ is $X$. If a different 
basepoint (or no 
basepoint at all)  is used, then Theorem 2 remains valid 
with possibly 
a smaller range of homotopy groups but Theorem 3 may fail 
completely. This is illustrated by the variety of Example 
1 above, where the 
results for all basepoints were analyzed in \cite{Gu3}.  
Nevertheless, it turns out 
that in all our counterexamples to (the conclusion of) 
Theorem 3, 
a modification of 
Theorem 3 holds in which $\Hol_D^\ast(S^2,X)$ is replaced 
by an appropriate 
subspace.  Full details will appear elsewhere.

The author is greatly indebted to Andrzej Kozlowski and 
Kohhei Yamaguchi for
their generous encouragement and collaboration.

\enddocument